\documentclass[12pt]{amsart}

\usepackage{amsthm}

\usepackage{graphicx}\usepackage{epsfig}
\usepackage{epstopdf}

\usepackage{latexsym}
\usepackage{amsmath}
\usepackage{amssymb}
\usepackage{amsfonts}

\theoremstyle{definition}

\newtheorem{dfn}{Definition}[section]
\newtheorem{defn}[dfn]{Definition}

\newtheorem{example}[dfn]{Example}
\newtheorem{rem}[dfn]{Remark}

\theoremstyle{plain}

\newtheorem{thm}[dfn]{Theorem}

\newtheorem{lem}[dfn]{Lemma}

\newtheorem{prop}[dfn]{Proposition}

\newtheorem{assumption}[dfn]{Assumption}
\newtheorem{cor}[dfn]{Corollary}

\newtheorem{conjecture}[dfn]{Conjecture}

\newtheorem{notation}[dfn]{Notation}

\def\proof{\par\medskip\noindent{\it Proof. }}

\def\R{{\mathbb R}}
\def\C{{\mathbb C}}
\def\Z{{\mathbb Z}}

\def\H{{\mathbb H}}

\def\O{{\mathcal O}}
\def\eps{\epsilon}
\def\al{\alpha}
\def\be{\beta}
\def\ga{\gamma}
\def\Ga{\Gamma}
\def\del{\delta}

\def\si{\sigma}

\def\la{\lambda}

\def\acts{\curvearrowright}

\def\embed{\hookrightarrow}
\def\8{\infty}
\def\<{\langle}
\def\>{\rangle}

\def\ol{\overline}

\def\t{\tilde}

\begin{document}

\title{Non-coherence of arithmetic hyperbolic lattices}
\author{Michael Kapovich}
\date{May 22, 2010}
\address{Department of Mathematics, 1 Shields Ave., University of California, Davis, CA 95618}
\email{kapovich@math.ucdavis.edu}

\begin{abstract}
We prove, under the assumption of the virtual fibration conjecture for arithmetic hyperbolic 3-manifolds, that all arithmetic lattices in $O(n,1)$, $n\ge 4, n\ne 7$, are non-coherent. We also establish noncoherence of uniform arithmetic lattices of the simplest type in $SU(n,1), n\ge 2$, and of  
uniform lattices in $SU(2,1)$ which have infinite abelianization. 
\end{abstract}

\maketitle

\section{Introduction}

Recall that a group $\Ga$ is called {\em coherent} if every finitely-generated subgroup of $\Ga$ is 
finitely-presented. This paper is motivated by the following 

\begin{conjecture}\label{conj:main}
Let $G$ be a semisimple Lie group which is not locally isomorphic to $SL(2,\R)$ and $SL(2, \C)$. 
Then every lattice in $G$ is noncoherent. 
\end{conjecture}

In the case of lattices in $O(n,1)$, this conjecture is due to Dani Wise. Conjecture \ref{conj:main} is true 
for all lattices containing direct product of two nonabelian free groups since the latter are incoherent. 
Therefore, it holds, for instance, for $SL(n, \Z), n\ge 4$. The case $n=3$, to the best of my knowledge, 
is unknown (this problem is due to Serre, see the list of problems \cite{Wall(1979)}).   

Conjecture \ref{conj:main} is out of reach for non-arithmetic lattices in $O(n,1)$ and $SU(n,1)$, since we do not understand the structure of such lattices. However, all known constructions of nonarithmetic 
lattices lead to noncoherent groups: See  \cite{KPV} for the case of Gromov--Piatetsky-Shapiro 
construction; the same argument proves noncoherence of 
non\-arith\-me\-tic reflection lattices (see e.g. \cite{Vinberg(1967)}) 
and non-arithmetic lattices obtained via Agol's \cite{Agol(2006)} construction. 
In the case of lattices in $PU(n,1)$, all known nonarithmetic groups are commensurable to the ones obtained via Deligne-Mostow construction \cite{Deligne-Mostow}. Such lattices contain fundamental groups of complex-hyperbolic surfaces which fiber over hyperbolic Riemann surfaces. Noncoherence of such groups is proven in \cite{Kapovich(1998a)}, see also section \ref{sec:ch}.  

In this paper we will discuss the case of arithmetic subgroups of rank 1 Lie groups. Conjecture \ref{conj:main} was proven in \cite{KPV} for non-uniform arithmetic lattices in $O(n,1), n\ge 6$ (namely, 
it was proven that the noncoherent examples from \cite{KP1} embed in such lattices). 
The proof of Conjecture \ref{conj:main}  in the case of all arithmetic lattices of the {\em simplest type} appears as a combination of \cite{KPV} and \cite{Agol(2008)}. In particular, it covers the case of all non-uniform arithmetic lattices ($n\ge 4$) and all 
arithmetic lattices in $O(n,1)$ for $n$ even, since they are of the simplest type. For odd $n\ne 3, 7$,   
there are also arithmetic lattices in $O(n,1)$ of ``quaternionic origin'' (see Section \ref{sec:quat} 
for the detailed definition), while for $n=7$ there is one more family of arithmetic groups associated with octonions. One of the keys to the proof of noncoherence above is virtual fibration theorem for various classes of hyperbolic 3-manifolds. Our main result (Theorem \ref{thm:main}) will be conditional to the existence of such fibrations:

\begin{assumption}\label{RFRS}
We will assume that every arithmetic hyperbolic 3-manifold $M$ of ``quaternionic origin'' admits a virtual fibration, i.e., 
$M$ has a finite cover which fibers over the circle. 
\end{assumption}

We discuss in Section \ref{vfc} what is currently known about the existence of virtual fibrations on such arithmetic manifolds. In short, 
assuming that a recent paper on subgroup separability by Dani Wise is correct, they all do.

%Note that 
%Ian Agol proved \cite{Agol(2008)} existence of such fibrations under the assumption that $\pi_1(M)$ satisfies virtual RFRS condition. 
%In  \cite{Agol(2009)} he proved RFRS condition for arithmetic hyperbolic 3-manifolds $M$, 
%whose fundamental groups are LERF. 

Our main result is

\begin{thm}\label{thm:main}
Under the assumption \ref{RFRS}, Conjecture \ref{conj:main} holds for all arithmetic lattices of quaternionic type. 
\end{thm}

In section \ref{sec:ch} we will also provide some partial corroboration to Conjecture \ref{conj:main}  for arithmetic subgroups of $SU(n,1)$. 
More precisely, we will prove it for uniform arithmetic lattices of the simplest type (also called {\em type 1 arithmetic lattices}) in $SU(n,1)$ 
and for all uniform lattices (arithmetic or not) in $SU(2,1)$ with (virtually) positive first Betti number. We will also prove the conjecture for all lattices in the isometry groups of the quaternionic--hyperbolic spaces ${\bf H}\H^n$ and the octonionic--hyperbolic plane ${\mathbf O}\H^2$. These groups are noncoherent since they contain arithmetic lattices from $O(4,1)$ or $O(8,1)$.

\medskip 
{\bf Historical Remarks.} The fact that the group $F_2\times F_2$ is non-coherent was known 
for a very long time (at least since \cite{Grunewald(1978)}). 
Moreover, it was proven by Baumslag and Roseblade \cite{Baumslag-Roseblade} that 
``most'' finitely generated subgroups of $F_2\times F_2$ are not finitely-presented. 
It therefore follows that many higher rank lattices (e.g., $SL(n, \Z)$, $n\ge 4$) are non-coherent. 
It was proven by P.~Scott \cite{Scott(1973a)} that finitely generated 3-manifold groups are all 
finitely-presented.  In particular, lattices in $SL(2,\C)$ are coherent. The first examples of incoherent geometrically finite groups in $SO(4,1)$ were constructed by the author and L.~Potyagailo \cite{KP1, Potyagailo(1990), Potyagailo(1994)}. These examples were generalized by B.~Bowditch and G.~Mess 
\cite{Bowditch-Mess} who constructed incoherent uniform arithmetic lattices in $SO(4,1)$. 
(For instance, the reflection group in the faces of right-angled 120-cell in $\H^4$ is one of such lattices.) Their examples, of course, embed in all other rank 1 Lie groups (except for $SO(n,1), SU(n,1), n=1, 2, 3$). Noncoherent arithmetic lattices in $SU(2,1)$ (and, hence, $SU(3,1)$) were constructed in \cite{Kapovich(1998a)}. All these constructions were ultimately based on either existence of hyperbolic 3-manifolds fibering over the circle (in the case of discrete subgroups in $SO(n,1)$) or existence of complex-hyperbolic surfaces which admit singular holomorphic fibrations 
over hyperbolic complex curves. A totally new source of 
noncoherent geometrically finite groups comes from the recent work of D.~Wise \cite{Wise(2008)}: He proved that fundamental groups of many (in some sense, most) polygons of finite groups are noncoherent. 
On the other hand, according to \cite{Kapovich(2005)}, the fundamental group of every even-sided (with at least 6 sides) hyperbolic polygons of finite groups embeds as a discrete convex-cocompact subgroup in some $O(n,1)$. 

{\bf Acknowledgments.} During the work on this paper I 
was supported by the NSF grants DMS-05-54349 and DMS-09-05802. I am grateful to Dani Wise for sharing with me the results \cite{Wise(2009)}, to Jonathan Hillmann for providing reference to Theorem \ref{hillman}, 
to Ben McReynolds and Andrei Rapinchuk for help with classification of arithmetic subgroups of $Isom({\mathbf H}\H^n)$ and 
$Isom({\mathbf O}\H^2)$ and to John Millson for providing me with references. I am also graterful to Leonid Potyagailo and Ernest Vinberg as this paper has  originated in our joint work \cite{KPV} on non-coherence of non-uniform arithmetic lattices.

\section{Relative separability for collections of subgroups}

Recall that a subgroup $H$ in a group $G$ is called {\em separable} if for every $g\in G\setminus H$ there exists a finite index subgroup $G'\subset G$ containing $H$ but not $g$. For instance, if $H=\{1\}$ then separability of $H$ amounts to residual finiteness of $G$. A group $G$ is called LERF if 
every finitely-generated subgroup of $G$ is separable. According to \cite{Scott(1978)}, LERF property is stable under group commensurability. 
Examples of LERF groups include:
Free groups (M.~Hall \cite{Hall(1949)}), surface groups (P.~Scott \cite{Scott(1978)}), certain hyperbolic 3-manifold groups (R.~Gitik \cite{Gitik(1999a)}, D.~Wise \cite{Wise(2006)}), groups commensurable to right-angled hyperbolic Coxeter  groups (P.~Scott \cite{Scott(1978)}, F.~Haglund \cite{Haglund(2008)}), 
non-uniform arithmetical lattices of small dimension (I.~Agol, D.~Long and A.~Reid \cite{Agol-Long-Reid} 
and M.~Kapovich, L.~Potyagailo, E.~B.~Vinberg \cite{KPV}), 
and other classes of groups (F.~Haglund and J.~{\'S}wi{\c{a}}tkowski \cite{Haglund-Swiat}, D.~Wise \cite{Wise(2000)}). Note that in a number of these results one has to replace finite generation of subgroups with geometric finiteness/quasiconvexity. On the other hand, there are 3-dimensional graph-manifolds whose fundamental groups are not LERF, see \cite{Long-Niblo}. For the purposes of this paper, we will need LERF property only  for the surface groups. 

\begin{rem}\label{rem:separability}
Dani Wise has a recent preprint \cite{Wise(2009)} which among other things proves that if $M$ is a closed 3-dimensional hyperbolic 
manifold containing an incompressible surface whose fundamental group is quasi-fuchsian, then $\pi_1(M)$ is LERF. His proof is rather long and intricate; whether or not it is correct is unclear at this point. For the purposes of our paper, we note the following. 
Suppose that $M$ is a closed 3-manifold, whose fundamental group is of quaternionic origin in the sense of section \ref{sec:quat}. Then, by applying \cite{Li-Millson} or \cite{Lubotzky(1996)}, see also \cite{Labesse-Schwermer} and the last paragraph of the introduction to 
\cite{Li-Millson}, one 
concludes that $M$ admits a finite cover $M'$ so that $b_1(M')\ge 2$ (actually, the first Betti number of a finite cover 
could be made arbitrarily large). It then follows from Thurston's results on Thurston's norm \cite{Thurston(1985)}, that $M'$  contains 
a nonseparating incompressible surface $\Sigma$ which is not a  fiber in a fibration over the circle. Therefore, $\pi_1(\Sigma)\subset \pi_1(M')\subset SL(2, \C)$ is a quasi-fuchsian subgroup. Hence, Wise's theorem applies to this class of arithmetic hyperbolic 3-manifolds. 
\end{rem}

 In this paper we will be using a {\em relative version} of subgroup separability, which deals with finite collections of subgroups of $G$. It is both stronger than subgroup separability (since it deals with collections of subgroups) and weaker than subgroup separability (since it does not require as much as separability in the case of a single subgroup). 
Actually, we will need this concept only in the case of pairs of subgroups, but we included the more general discussion for the sake of completeness.

Let $G$ be a group, $H_1, H_2\subset G$ be subgroups. We say that a double coset 
$H_1 g H_2$ is {\em trivial} if it equals the double coset $H_1 \cdot 1\cdot  H_2=H_1 \cdot   H_2$. 
In other words, $g\in H_1\cup H_2$. 
Given a finitely-generated group $G$, we let $\Ga_G$ denote its Cayley graph (here we are abusing the notation by suppressing the choice of a finite generating set which will be irrelevant for our purposes). Recall also 
that a {\em geometric action} of a group $G$ on a metric space $X$ is an isometric properly discontinuous cocompact action. 

\begin{defn}
Let $G$ be a group, ${\mathcal H}=\{H_1, H_2, ...H_M\}$ be a collection of subgroups of $G$. We will say that ${\mathcal H}$ is {\em relatively separable} in $G$ if for every finite collection of nontrivial double cosets $H_i g_k H_{j}, i, j\in \{1, ...,M\}, k=1,...,K$,  there exists a finite index subgroup $G'\subset G$ which is disjoint from the above double cosets.  
\end{defn}

In the case when $G, H_i, i=1,...,M$, are finitely generated, separability can be reformulated as follows: 

Given a number $R$, there exists a finite index subgroup $G'\subset G$ so that for each $g\in G'$ 
either $g\in H_i \cdot H_j$ for some $i, j$ (and, hence, $d(\Ga_{H_i}, g\Ga_{H_j})= d(\Ga_{H_i}, \Ga_{H_j})$) or 
$$
d(g \Ga_{H_i}, \Ga_{H_j})\ge R
$$
for all $i, j\in \{1,....,M\}$. 

Equivalently, in the above property one can replace $\Ga_{G}$ with a space $X$ on which $G$ acts geometrically and $\Ga_{H_i}$'s with $H_i$-invariant subsets $X_i\subset X$ with compact quotients 
$X_i/H_i$, $i=1,...,M$. 

In what follows we will use the following notation

\begin{notation}
Given a metric space $X$, a subset $Y\subset X$ and a real number $R\ge 0$, let $B_R(Y):=\{x\in X: d(x, Y)\le R\}$, i.e., the $R$-neighborhood of $Y$. For instance for $x\in X$, $B_R(x)$ is the closed $R$-ball in $X$ centered at $x$.  We let 
$proj_Y: X\to Y$ denote the nearest-point projection of $X$ to $Y$. 
\end{notation}

 Below are several useful examples illustrating relative separability.

\begin{example}
Suppose that $M=1$. Then ${\mathcal H}$ is relatively separable provided that $H=H_1$ is separable in 
$G$.  
\end{example}
\proof Suppose that $H$ is separable in $G$. Given $R$, there  are only finitely many distinct nontrivial double cosets 
$$
H g_k H, k=1,...,K,
$$
so that 
$$
d(g \Ga_{H}, \Ga_{H})< R
$$
for $g\in H g_k H$. Let $G'\subset G$ be a subgroup containing $H$ but not $g_1,...,g_K$. Then 
$G'$ is disjoint from 
$$
\bigcup_{k=1}^K H g_k H
$$
and the claim follows. \qed

Although the converse to the above example is, probably, false, relative separability suffices for typical applications of subgroup separability. For instance, suppose that $G$ is the fundamental group of a 3-manifold $M$ and ${\mathcal H}=\{H\}$ is a relatively separable surface  subgroup of $G$. Then 
a finite cover of $M$ contains an incompressible surface (whose fundamental group is a finite index subgroup in $H$).

\begin{example}
Let $G$ be an arithmetic lattice of the simplest type in $O(n,1)$ and $H_1,..., H_M\subset G$ be 
the stabilizers of distinct ``rational'' hyperplanes $L_1,...,L_M$ in $\H^n$, i.e., $L_i/H_i$ has finite volume, $i=1,...,M$. Then  ${\mathcal H}=\{H_1,..., H_M\}$ is relatively separable in $G$. See \cite{KPV}. 
\end{example}

\begin{example}
As a special case of the second example, suppose that $G$ is a surface group and $H_1,..., H_M$ are  cyclic subgroups. Then  ${\mathcal H}=\{H_1,..., H_M\}$ is relatively separable in $G$. 
\end{example}

Recall that finitely generated subgroups of surface groups are separable. Therefore, 
one can generalize the last example as follows:

\begin{prop}\label{prop:sep}
Suppose that $G$ is a word-hyperbolic group which is separable with respect to its quasiconvex subgroups. Let $H_1,..., H_M$ be residually finite quasiconvex subgroups with finite pairwise intersections. Then ${\mathcal H}=\{H_1,..., H_M\}$ is relatively separable in $G$. 
\end{prop}
\proof Let $H\subset G$ be a subgroup generated by sufficiently deep finite index torsion-free 
subgroups $H_i'\subset H_i$ ($i=1,...,M$). Then, according to \cite{Gitik(1999b)}, $H$ is isomorphic to the free 
product $H_1'*H_2' *... * H_M'$  and is quasiconvex. Let $X$ denote the Cayley graph of $G$ and $X_i, i=1,...,M$ the Cayley graphs of $H_1,...,H_M$. Since the groups $H_i$ are quasiconvex and have finite intersections, for every $R<\infty$, there exists $r$ so that for $i\ne j$, 
the projection of $B_R(X_i)$ to $X_j$ is contained in $B_{r}(1)$. 

Let $Y_j$ denote the preimage in $X$ of $B_{r}(1)$ under the nearest-point projection $X\to X_j$. Thus, 
$B_R(X_j)\subset Y_i, i\ne j$. Moreover, if we choose $r$ large enough then 
$$
Y_j^c\cap Y_i^c=\emptyset, \quad \forall i\ne j,$$
where $Y_j^c$ denotes the complement of $Y_j$ in $X$. 

Since the group $H_j$ is residually finite, there exists a finite index subgroup $H_j'\subset H_j$ so that 
$Y_j$ is a {\em sub-fundamental domain} for the action $H_j'\acts X$, i.e.,
$$
h(Y_j)\cap Y_j=\emptyset, \forall h\in H_j'\setminus \{1\}. 
$$

Therefore, one can apply the ping-pong arguments to the collection of 
sub-funda\-men\-tal domains $Y_1,...,Y_M$ as follows: 

Every nontrivial element $h$ is the product 
$$
h_{i_1}\circ h_{i_2}\circ ... \circ h_{i_m}
$$
where $h_{i_k}\in H_{i_k}\setminus \{1\}$  
and $i_k\ne i_{k+1}$ for each $k=1,...,m-1$. Then, arguing inductively on $m$, 
we see that for each $j$, 
$$
h(Y_j)\subset Y_{l}^c,  
$$
where $l=i_1$. 

We now claim that $d(X_i, h(X_j))\ge R$ provided that $h\notin H_i \cdot H_j$. 
We write down $h$ in the normal form as above with $l=i_1$. 

1. Suppose first that $i\ne l$. Then $B_R(X_i)\subset Y_l$. On the other hand, by taking any $m\ne j$, we get 
$$
h(X_j) \subset h(Y_m)\subset Y_{l}^c. 
$$
Thus $B_R(X_i)\subset X_l$ has empty intersection with $h(X_j)\subset Y_l^c$ and the claim follows. 

2. Suppose now that $i=l=i_1$. Then $s=i_2\ne i$; set
$$
g:= h_{i_2}\circ ... \circ h_{i_m}.
$$
Then 
$$
d(X_i, h(X_j))= d(X_i, h_i g (X_j))= d(X_i, g(X_j)). 
$$
Now, by appealing to Case 1, we get
$$
d(X_i, h(X_j))=d(X_i, g(X_j))\ge R. 
$$
We hence conclude that ${\mathcal H}$ is relatively separable in the subgroup generated by 
$H_1,....,H_M$. 

There are only finitely many nontrivial double coset classes $H_i g_k H_j$, $k=1,...,K,$ 
in $H_i\backslash G /H_j$ so that for the elements $g\in H_i g_k H_j\subset G$, we have
$$
d(X_i, g(X_j))< R.
$$
Note that $g_k\notin H$ unless $g_k\in H_i \cdot H_j$ (in which case the corresponding 
double coset would be trivial). Since $H$ is quasiconvex in $G$ and does not contain $g_k, k=1,...,K$, 
by the subgroup separability of $G$, there exists a finite index subgroup $G'\subset G$ containing $H$, so that $g_1,...,g_K\notin G'$. Therefore, $G'$ has empty intersection with each of the double cosets 
$H_i g_k H_j, k=1,...,K$ (for all $i, j\in \{1,...,M\}$). 
It therefore follows that for every $g\in G'\setminus H_i \cdot H_j$, 
$$
d(X_i, g(X_j))\ge R. 
$$
Hence ${\mathcal H}$ is relatively  separable in $G$. \qed 

\medskip
In order to apply this proposition in section \ref{final}, we need the following simple

\begin{lem}\label{easy}
Let $N\triangleleft G$ be a normal finitely-generated subgroup and ${\mathcal H}$ be a finite set of subgroups 
in $N$ which is relatively separable in $N$. Assume that $G/N$ is residually finite. 
Then  ${\mathcal H}$ s relatively separable in $G$.
\end{lem}
\proof Let $H_i g_k H_j$, $g_k\in A\subset G$, 
be the double cosets that we want to avoid. We will construct 
the appropriate finite-index subgroup in $G$ in two steps. First, consider those 
double cosets $H_i g_k H_j, g_k\in A_1$, which are contained in $N$. By separability of ${\mathcal H}$ in $N$, there exists a finite-index subgroup $N'\subset N$ which is disjoint from these double cosets. This subgroup may not be normal in $G$, but 
it is standard that there exists a finite-index subgroup $G_1\subset G$ containing $N_1$ (the subgroup $G_1$ is the normalizer of $N_1$ in $G$). 
Hence, $G_1$ is disjoint from the first set of double cosets. Now, consider double cosets 
$H_i g_k H_j, g_k\in A_2$, so that $g_k\notin N$. Then, since $G/N$ is residually finite, there exists 
a finite-index subgroup $G_2\subset G$ so that $g_k\notin G_2, g_k\in A_2$, but $N\subset G_2$. 
Clearly, $G_2$ is disjoint from the double cosets $H_i g_k H_j, g_k\in A_2$. Now, taking 
$G'=G_1\cap G_2$ we obtain a finite-index subgroup in $G$ which is disjoint from 
the double cosets $H_i g_k H_j, g_k\in A_1\cup A_2$.\qed

\section{Virtual fibration conjecture}\label{vfc}

W.~Thurston conjectured that every closed hyperbolic 3-manifold $M$ is {\em virtually fibered}, i.e., 
it admits a finite cover $M'$ which is fibered over the circle. This conjecture, known as {\em Thurston's Virtual Fibration Conjecture} 
was proven by I.~Agol in \cite{Agol(2008)} under the assumption that $\pi_1(M)$ is virtually RFRS. Moreover, under this 
assumption, $M$ admits infinitely many non-isotopic virtual fibrations. 

The latter condition holds provided that $\pi_1(M)$ contains a finite-index subgroup which embeds in a right-angled  Coxeter (or Artin) group. This allowed Agol to prove \cite{Agol(2008)} the virtual fibration conjecture for all arithmetic manifolds of the simplest type: Immersed closed totally-geodesic surfaces provide a way to virtually embed $\pi_1(M)$ in a right-angled Artin group using the results of Haglund and Wise \cite{Haglund-Wise}. 
Paper \cite{Wise(2009)} of Dani Wise claims virtual RFRS condition for fundamental groups of all closed hyperbolic 3-manifolds containing 
an incompressible surface with quasi-fuchsian fundamental group. Recall that every  closed arithmetic 3-manifold of quaternionic origin admits a finite cover 
containing such an incompressible surface, see Remark \ref{rem:separability}. 
Thus, assuming Dani Wise's results \cite{Wise(2009)}, if $M$ is a closed arithmetic 3-manifold of quaternionic origin, then 
Virtual Fibration Conjecture holds for $M$. 

We also note that if $M$ is arithmetic, then existence of one virtual fibration implies existence of infinitely many virtual fibrations since the commensurator $Comm(\Ga)$ of $\Ga=\pi_1(M)$ is dense in $PSL(2,\C)$: If $F\triangleleft \Ga$ is a normal surface subgroup, there exists $\al\in Comm(\Ga)$ so that $F':=\al F \al^{-1}$ is not commensurable to $F$. Therefore, $F'\cap \Ga$ is a surface subgroup corresponding 
to a different fibration of the manifold $M'=\H^3/\Ga'$, where $\Ga':= \Ga \cap \al \Ga \al^{-1}$.

\section{Normal subgroups of Poincar\'e duality groups}\label{Appendix}

Recall that a group $\Gamma$ is called an $n$-dimensional 
{\em Poincar\'e Duality group} (over $\Z$), abbreviated $PD(n)$ group, if 
there exists $z\in H_n(\Gamma, D)$, so that 
$$
\cap z: H^i(\Gamma, M)\to H_{n-i}(\Gamma, \overline{M})
$$ 
is an isomorphism for $i=0,..,n$ and every $\Z\Ga$--module $M$. Here   
$\overline{M}= D\otimes M$, where $D\cong H^n(\Gamma, \Z\Gamma)$ is 
the {\em dualizing module}. For instance, if $X$ is a closed $n$-manifold so that 
$X=K(\Ga,1)$, then $\Ga$ is a $PD(n)$ group. The converse holds for $n=2$, while for $n\ge 3$ 
the converse is an important open problem (for groups which admit finite $K(\Gamma,1)$). 

Recall also that a group $\Ga$ is called $FP_r$ (over $\Z$) if there exists a partial resolution 
$$
P_r \to P_{r-1}\to ... \to P_0 \to \Z \to 0
$$ 
of finitely-generated projective $\Z\Gamma$-modules. For instance, $\Gamma$ is $FP_1$ iff it is 
finitely--generated, while every finitely--presented group is $FP_2$ (the converse is false, see \cite{Bestvina-Brady}). We refer the reader to \cite{Bieri(1976a)} for a comprehensive discussion of 
$PD(n)$ and $FP_r$ groups. We will need the following theorem in the case $n=4$, $r=2$: 

\begin{thm}\label{hillman}
(J.~Hillman \cite[Theorem 1.19]{Hillman(FMGK)}, see also \cite{Hillman-Kochloukova}.) 
Let $\pi$ be a $PD(n)$-group with an $FP_r$ normal subgroup $K$ 
such that $G=\pi/K$ is a $PD(n-r)$ group and $2r\geq n-1$.
Then $K$ is a $PD(r)$-group.
\end{thm}

\section{Bisectors}\label{bisectors}

Given $p\ne q\in \H^n$ the {\em bisector} $Bis(p,q)$ is the set of points in $\H^n$ 
equidistant from $p$ and $q$. Define the closed half-spaces $Bis(p,q)^\pm$ bounded by $Bis(p,q)$ by requiring $Bis(p,q)^+$ to contain  $q$ and $Bis(p,q)^-$ to contain $p$.

Consider a configuration $C$ of 3-dimensional subspaces  $H_1, H_2, H_3$  in 
$\H^n$ and geodesics $\ga_1,...,\ga_4$, so that:

1. $\ga_i\subset H_i, i=1, 2, 3, \ga_4\subset H_3$.

2. $H_1\cap H_2=\ga_2, H_2\cap H_3=\ga_3$ and the angles at 
these intersections are $\ge \al>0$. 

3. $\ga_i, \ga_{i+1}$ ($i=1, 2, 3$) are 
at least distance $R_0>0$ apart. 

The following is elementary: 

\begin{lem}\label{L*}
Given $\al, r$ and $R_0$, there exists $R_*$ so that if $d(\ga_2, \ga_{3})\ge R_*$ then 
$d(H_1, H_3)\ge r$. 
\end{lem}

Therefore, from now on, we will also fix a number $r>0$ so that the configuration $C$ satisfies: 

4. $H_1$, $H_3$ are at least distance $r>0$ apart.

\medskip 
Let $p_i, q_i\in \ga_i$ denote the points 
closest to $\ga_{i-1}, \ga_{i+1}$ respectively. 
We let $m_i$ denote the midpoint of $\ol{q_i p_{i+1}}, i=1, 3$. 
Let $n_i$ denote the midpoint of $\ol{p_i q_i}, i=2, 3$. 
Our assumptions imply that  $d(q_i, p_{i+1})\ge R_0, i=1,...,3$.  

\begin{figure}[htbp]
\begin{center}
 \includegraphics[width=3.5in]{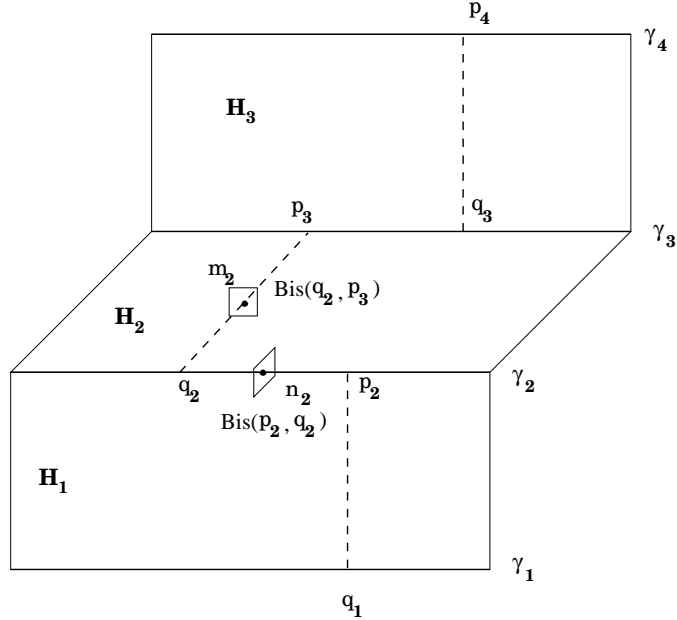}
\caption{{\bf Configuration of planes and lines.}}
\label{Fig1}
\end{center}
\end{figure}

See Figure \ref{Fig1}. 

\begin{lem}\label{L1}
There exist $R_1=R_1(\al, r, R_0), \rho_1=\rho_1(\al, r, R_0)>0$ so that if $d(q_2, p_2)\ge R_1$ then: 

1.  The bisectors $Bis(p_2, q_2)$ and $Bis(q_2, p_3)$ are within distance. 

2.  
$$
\ga_1\subset Bis(p_2, q_2)^-, H_3\subset Bis(p_2, q_2)^+. 
$$

3. $d(\ga_1, Bis(p_2, q_2))\ge \rho_1$ and $d(H_3, Bis(p_2, q_2))\ge \rho_1$. 
%In particular, the bisectors  $Bis(p_2, q_2)$ and $Bis(q_2, p_3)$ 
%intersect the path $$\eta= n_2 q_2 m_2, $$ only at the end-points of $\eta$. 
See Figure \ref{Fig1}. 
\end{lem}
\proof Observe that as $d=d(q_2, p_2)$ goes to infinity, the bisector 
$Bis(p_2, q_2)$  converges to an ideal point of the geodesic $\ga_2$. 
Since both $\ga_2, Bis(q_2, p_3)$  are orthogonal to $\ol{q_2 p_3}$ 
and $d(p_3, q_2)\ge R_0>0$,  it follows that 
$Bis(p_2, q_2), Bis(q_2, p_3)$ are disjoint for large $d$. Similarly, 
since $\ga_2$ is within positive distance from $H_3$, it follows that 
$H_3\subset Bis(p_2, q_2)^+$ for large $d$. Similarly, $\ga_1\subset Bis(p_2, q_2)^-$ 
for large $d$ (since $d(q_1, p_2)\ge R_0>0$). Alternatively, one compute $R_1, \rho_1$ directly 
using hyperbolic trigonometry. \qed

The next two lemmas are proven by the same arguments as Lemma \ref{L1}.

\begin{lem}\label{L2}
Assume that $d(q_2, p_3)= R_0$. Then there exist $R_2=R_2(\al, r, R_0), 
\rho_2=\rho_2(\al, r, R_0)>0$ so that 
if $d(p_2,q_2)\ge R_1, d(p_3, q_3)< R_1$ (where $R_1$ is as in Lemma \ref{L1}), 
and $d(q_3, p_4)\ge R_2$, then:

1. The bisectors $Bis(p_2, q_2)$ and $Bis(q_3, p_4)$ are within positive distance.

2.  $H_1\cup H_2\subset Bis(q_3, p_4)^-$. 

3. $d(H_1\cup H_2, Bis(q_3, p_4))\ge \rho_2$. 

%In particular, the bisectors  $Bis(p_2, q_2)$ and $Bis(q_3, p_4)$  
%intersect the path
%$$\eta= n_2 q_2 p_3 q_3 m_3, $$only at the end-points of $\eta$. 
\end{lem}

\begin{lem}\label{L3}
Assume again that $d(q_2, p_3)= R_0$ and, in addition, 
$d(p_2,q_2)< R_1$, $d(p_3, q_3)< R_1$. 
Then there exist $R_3=R_3(\al, r, R_0), \rho_3=\rho_3(\al, r, R_0)>0$ so that if 
$$
d(q_1, p_2)\ge R_3, \quad d(q_3, p_4)\ge R_4,
$$
then:

1. The bisectors $Bis(q_1, p_2)$ and $Bis(q_3, p_4)$ are within distance. % $\ge \rho_3$. 

2. $H_2\cup H_3\subset Bis(q_1, p_2)^+$ and $H_2\cup H_1\subset Bis(q_3, p_4)^-$. 

3. $d(H_2\cup H_3, Bis(q_1, p_2))\ge \rho_3$, 
$d(H_2\cup H_1, Bis(q_3, p_4))\ge \rho_3$.

%In particular, 
%the bisectors  $Bis(q_1, p_2)$ and $Bis(q_3, p_4)$ 
%intersect the path
%$$\eta= m_1 p_2 q_2 p_3 q_3 m_3, $$
%only at the end-points of $\eta$. 
\end{lem}

We now set $R_4:=\max(R_2, R_3)$ and $\rho_4:=\min(\rho_1, \rho_2, \rho_3)$. 
Thus, $R_4=R_4(\al, r, R_0)$  and $\rho_4:=\rho_4(\al, r, R_0)$. 

We will say that $H_i$ is {\em large}  (relative to the configuration $C$) if $d(q_i, p_{i+1})\ge R_4$ 
and $\ga_i$ is {\em large} if   $d(p_i, q_{i})\ge R_1$. If $H_i$ or $\ga_i$ is not large, we will call 
it {\em small}.

\section{A combination theorem for quadrilaterals of groups}\label{sec:comb}

Combination theorems in theory of Kleinian groups provide a tool for proving that a 
subgroup of $O(n,1)$ generated by certain discrete subgroups is also discrete and, moreover, has 
prescribed algebraic structure. The earliest example of such theorem is ``Schottky construction'' 
(actually, due to F.~Klein) producing free discrete subgroups of $O(n,1)$. This was generalized by B.~Maskit in the form of Klein-Maskit combination theorems where one constructs amalgamated free products and HNN extensions acting properly discontinuously on $\H^n$. More generally, the same line of arguments applies to graphs of groups. Complexes of groups are higher-dimensional generalizations of graphs of groups. A combination theorem for polygons of finite groups was proven in \cite{Kapovich(2005)}. The goal of this section is to prove a combination theorem for certain quadrilaterals of infinite groups.

Let $G_1$ be a discrete subgroup in $Isom(L)$, $L\cong \H^3$. Pick two nonconjugate maximal cyclic subgroups $G_{\eps_1}, G_{\eps_2}$ in $G_1$. For $i=1, 2$, let $\ga_i=L(\eps_i)\subset L$ denote the  
invariant geodesics  of $G_{\eps_i}$.

We will assume that:  

\begin{assumption}\label{as1}
1. The distance between these 
geodesics is $R_0>0$. 

2. There exists $R>R_0$ such that for all distinct 
geodesics $\be, \ga$ in the $G_1$-orbits of $\ga_1, \ga_2$, the distance $d(\be, \ga)$ is at least $R$ 
unless there exists $g\in G_1$ which carries $\be\cup \ga$ to $\ga_1\cup \ga_2$. 

3. There exists $R_1>0$ so that for each $\ga=\ga_i, i=1, 2$ and geodesics 
$\be_1, \be_2\in G_0\cdot (\ga_1 \cup \ga_2)$ 
so that $d(\be_j, \ga)= R_0$ ($j=1, 2$), it follows that 
$$
d(proj_{\ga}(\be_1), proj_{\ga}(\be_2))\ge R_1.
$$

We will specify the choices of $R$ and $R_1$ later on. 
\end{assumption}

\begin{figure}[htbp]
\begin{center}
 \includegraphics[width=3.5in]{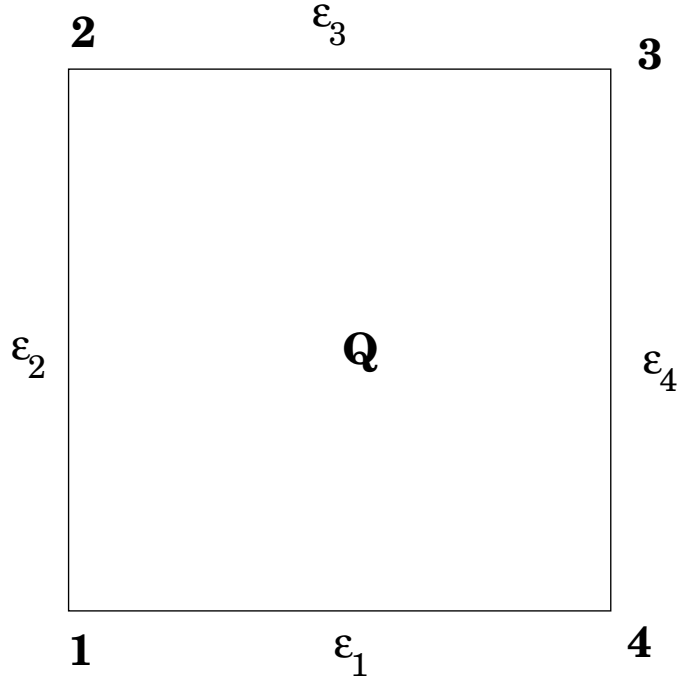}
\caption{{\bf Square of groups.}}
\label{Fig2}
\end{center}
\end{figure}

We then define a quadrilateral $Q$ of groups where the vertex groups $G_1,...,G_4$ are copies of 
$G_1$ and the edge groups $G_{\eps_1},...,G_{\eps_4}$ are copies of $G_{\eps_1}, G_{\eps_2}$. In particular, $Q$ has trivial face-group. We refer the reader to \cite{Bridson-Haefliger} for the precise definitions of complexes of groups and their fundamental groups, we note here only that for each vertex $v$ of $Q$ incident to an edge $e$, the structure of a quadrilateral of groups prescribes an embedding $G_e\embed G_v$.  The fundamental group $G=\pi_1(Q)$ of this quadrilateral of groups is the direct limit of the diagram of groups and homomorphisms given by $Q$.

More specifically, we require that $Q$ admits a $\Z_2\times \Z_2$--action, generated by two involutions $\si_1, \si_2$, so that
$$
\si_2(G_1)=G_2, \si_1(G_1)=G_4,
$$
and $\si_i$ fixes $G_{\eps_i}$, $i=1, 2$. The isomorphisms 
$$
G_1\to G_i, \quad i=2,3,4$$
 are induced by $\si_1, \si_1\si_2, \si_2$ respectively.

Let $X$ denote the universal cover of the complex of groups $Q$. Then $X$ is a square complex where the links of vertices are bipartite graphs. Since $G_{\eps_1}\cap G_{\eps_2} =\{1\}$, it follows that the links of $X$ contain no bigons. Thus, $X$ is a CAT(0) square complex.  Let $X^1$ denote 
1-skeleton of the complex $X$; we metrize  $X^1$ by declaring every edge to have unit lengths. 
Recall that for $k>0$ a path $\eta$ in $X^1$ is a $k$-local geodesic if 
every sub-path of length $k$ in $\eta$ is a geodesic in $X^1$.
%We consider the We identify paths in $X^1$ with the corresponding sequences of vertices. 

Note that if we choose edge groups which generate a free subgroup of $G$, then the links of the vertices of $X$ are trees.  Define the group $\tilde{G}=\<G, \si_1, \si_2\>$ generated by $G$, $\si_1, \si_2$. Then $\t{G}$ is a finite extension of $G$. The group $\t{G}$ acts on $X$ with exactly two orbit types of edges. For every edge $e$ we define $Type(e)=\eps_i$ if 
$g(e)$ projects (under the map $X\to Q$) to $\eps_i$ for some  $g\in \t{G}$. 

\medskip 
We next describe a construction of representations $\phi: G\to O(5,1)$ which are discrete and faithful 
provided that $R$ is sufficiently large and $R_1$ is chosen appropriately. 
Let $S_i$ denote 3-dimensional subspaces in $\H^5$ 
which intersect $L$ along $\ga_i$ at the angles $\al_i\ge \al>0$, $i=1, 2$. We assume that $S_1\cap S_2$ is a geodesic which is within distance $\ge \frac{r}{2}>0$ 
from $L$ and that $S_1$ is orthogonal to $S_2$. 
Let $\si_i, i=1,2$ denote commuting isometric involutions in $\H^5$ with the fixed-point sets $S_i, i=1,2$ 
respectively. 

Let $L_1:=L, L_2:=\si_2(L_1), L_3:=\si_1 \si_2(L_1), L_4:= \si_1(L_1)$. Then 

\begin{equation}\label{dist=r}
d(L_1, L_3)=d(L_2, L_4)\ge r. 
\end{equation}

We have the (identity) discrete embedding $\phi_1: G_1\to Isom(L_1)\subset Isom(\H^5)$. 
We will assume that $\phi_1(G_{\eps_i})$ stabilizes $\ga_i, i=1, 2$. 

Given this data, we define a representation $\phi: G\to Isom(\H^5)$ so that the symmetries 
$\si_1, \si_2$ of the quadrilateral of groups $Q$ correspond to the involutions $\si_1, \si_2\in O(5,1)$:  

1. $\phi_1=\phi|G_1$. 

2. $\phi_2=\rho|G_2= Ad_{\si_2}\circ \phi_1$,  $\phi_4=\phi|G_4= Ad_{\si_1}\circ \phi_1$. 

3. $\phi|G_3= Ad_{\si_1\si_2} \circ \phi_1$. 

Thus $\phi$ extends to a representation (also called $\phi$) of the group $\tilde{G}=\<G, \si_1, \si_2\>$ 
generated by $G$, $\si_1, \si_2$.

\begin{figure}[htbp]
\begin{center}
 \includegraphics[width=3.5in]{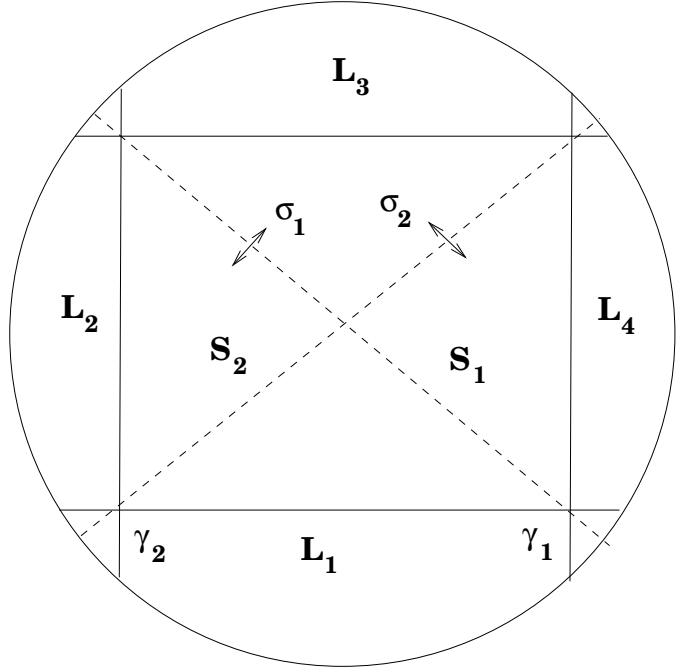}
\caption{{\bf Constructing a representation $\phi$: Projective model of $\H^5$.}}
\label{Fig3}
\end{center}
\end{figure}

Our main result is

\begin{thm}\label{comb}
If in the above construction $R, R_1$ is sufficiently large (with fixed $R_0$, $r$, $\al$), 
then $\phi$ is discrete and faithful. 
\end{thm}
\proof Our proof is analogous to the one in \cite{Kapovich(2005)}. 
Every vertex $x$ of $X$ is associated with a 3-dimensional hyperbolic subspace $L(x)\subset \H^5$, namely, it is a subspace stabilized by the vertex group of $x$ in $G$. If $x=g(x_1), g\in \t{G}$, 
where $x_1\in Q$ is stabilized by $G_1$, then $L(x)=g(L_1)$. 

Similarly, every edge $e$ of $X$ is associated a geodesic $L(e)\subset \H^5$ stabilized by 
$G_e$. Hence, if $e$ connects vertices $x, y$ then
$$
L(e)=L(x)\cap L(y). 
$$

The following properties of the subspaces $L(x), L(e)$ follow directly from the Assumption \ref{as1} 
and inequality \eqref{dist=r}:
 
If $x, y$ are vertices of $X$ which belong to a common 2-face and are 
within distance 2 in $X^1$, the subspaces $L(x), L(y)$ are within positive distance 
$r$ from each other. If $e, f$ are distinct edges incident to a common vertex $x$ of $X$ then 

a) either $e, f$ belong to a common 2-face of $X$, in which case 
$$
d(L(e), L(f))=R_0,$$

b) otherwise, $d(L(e), L(f))\ge R$. 

\medskip
In order to prove Theorem \ref{comb}, it suffices to show that for sufficiently large $R$, 
there exists $\del>0$ so that for all distinct edges $f, f'$ in $X$, 
\begin{equation}\label{rho}
d(L(f), L(f'))\ge \del. 
\end{equation}

Consider a {\em 4-chain} ${\mathcal C}$ of edges $(f_1,...,f_4)$ in $X^1$ so that 
$f_i\cap f_{i+1}=y_i, i=1, 2, 3$, and the concatenation $f_1\cup f_2 \cup f_3 \cup f_4$ is a local 
3-geodesic in $X^1$. Two chains $(f_1,...,f_4)$ and $(f_1',...,f_4')$ 
 are said to be {\em combinatorially isomorphic} if:  
 
1. $Type(f_i)=Type(f_i'), i=1,...,4$, and  

2. $f_i, f_{i+1}$ belong to a common 2-face iff $f_i, f_{i+1}$ belong to a common 2-face for 
$i=1,2,3$. 

\medskip 
Clearly, there are only finitely many combinatorial isomorphism types of 4-chains in $X^1$. 

Each 4-chain ${\mathcal C}=(f_1,...,f_4)$ corresponds to a configuration $C$ of three hyperbolic 3-dimensional  subspaces $L(y_1), L(y_2), L(y_3)$, and four 
geodesics $L(f_1),...,L(f_4)$ as in section \ref{bisectors}. The isometry class of $C$ depends only 
on the combinatorial isomorphism class of  ${\mathcal C}$. 

We now describe the choices of $R_1$ and $R$:

1. We will assume that $R_1$ is such  that 
for each combinatorial isomorphism class of 4-chains in $X^1$, the corresponding configuration $C$ 
of hyperbolic subspaces in $\H^5$ satisfies $R_1\ge R_1=R_1(C)$ where the function $R_1$ is defined 
in section \ref{bisectors}.

2. We will assume that $R$ is chosen so that 
for each combinatorial isomorphism class of 4-chains in $X^1$, the corresponding configuration $C$ 
of hyperbolic subspaces in $\H^5$ satisfies $R\ge R_4=R_4(C)$; we let  
$\rho$ be the minimum of the numbers $\rho_4(C)$, where the minimum is 
taken over all combinatorial isomorphism classes of 4-chains in $X^1$. Here $R_4$ and $\rho_4$ 
are the functions introduced in section \ref{bisectors}. 

\medskip 
We claim that one can take $\del=\min(2\rho, r, R_0)$ in \eqref{rho}. 

If $f, f'$ share a vertex then 
$d(L(f), L(f'))\ge R_0$. If $e, f$ do not share a vertex and belong to a common 2-face then 
$d(L(f), L(f'))\ge r$. We therefore consider the generic case when none of the above occurs. 

Consider a geodesic path $\eta$ in $X^1$ connecting $f$ to $f'$; 
the path $\eta$ is a concatenation of the edges 
$$
\eta=  e_2 \cup... \cup e_{k-1}.
$$
 Then $k\ge 3$ and for $e_1=f, e_{k}=f'$, the concatenation
 $$
 \theta= e_1 \cup \eta \cup e_{k}
 $$
is a 3-local geodesic.  
Set $e_{i}=[x_{i} x_{i+1}]$, $i=1,...,k$. Following section \ref{bisectors}, for each geodesic $L(e_i)$ 
we define points $p_i$ (closest to $L_{e_{i-1}}$) and $q_i$  (closest to $L_{e_{i+1}}$). 

We will say that an edge $e$ in $\theta$ is {\em large} if the corresponding geodesic $L(e)$ is large 
in the sense of section \ref{bisectors}. Similarly, a vertex $x_i$ in $\theta$ is {\em large} if 
$L(x_i)$ is large in the sense of section \ref{bisectors}, i.e., $d(q_i, p_{i+1})\ge R$. Accordingly, $x_i$ is small if $d(q_i, p_{i+1})=R_0$ and $e$ is small if it is not large. 

The next lemma immediately follows from Assumption \ref{as1}:

\begin{lem}\label{large exists}
For every edge $e_i=[x_{i} , x_{i+1}]$, at least one of $e_i, x_i, x_{i+1}$ is large. 
\end{lem}

\begin{rem}
Note that we can have a vertex $x_i$ so that all three $x_i, e_i, e_{i-1}$ are small.  
\end{rem}

We also define a sequence of bisectors in $\H^n$ corresponding to the path $\theta$ as follows:

For each large cell $c=e_i, c=x_i$ in $\theta$ we take the corresponding bisector $Bis(c)=Bis(p_i, q_i), Bis(c)=Bis(q_i, p_{i+1})$ defined as in section \ref{bisectors}. Recall that we also have half-spaces 
$Bis(c)^\pm\subset \H^n$ bounded by the bisectors $Bis(c)$.

We define the natural total order $>$ on the cells in $\theta$ by requiring that $e_j> e_i, x_j> x_i$ 
if $j> i$ and that $e_i> x_i$. 

\begin{prop}\label{bisector chain}
%1.  $Bis(x_i)$ separates $L(e_{i-1})$ from $L(e_{i})$, while $Bis(e_i)$ separates $L(e_{i-1})$ from 
%$L(e_{i+1})$. 

1.  If $c>c'$ then $Bis(c')^+ \subset Bis(c)^+$. 

2. If $e_i\ne c$ then $d(L(e_i), Bis(c))\ge \rho>0$. Moreover, if $e_i<c$ then 
$L(e_i)\subset Bis(c)^-$, while if $e_i>c$ then $L(e_i)\subset Bis(c)^+$. 
\end{prop}
\proof %1. The first assertion follows immediately from Lemmas \ref{L1}, \ref{L2}, \ref{L3}. 
Part 1. If $c, c'$ are consecutive (with respect to the order $>$) large cells in $\theta$ then they 
belong to a common 4-chain. Therefore, the assertion follows from Lemmas \ref{L1}, \ref{L2}, \ref{L3} 
in this case. The general case follows by induction since for three consecutive large cells 
$c_1< c_2< c_3$ in $\theta$, we have
$$
Bis(c_3)^+ \subset Bis(c_2)^+ \subset Bis(c_1)^+. 
$$

Part 2. (a) If $e_i, c$ are no separated (with respect to the order $<$) by any large cells, then 
they belong to a common 4-chain and the assertion follows from  
Lemmas \ref{L1}, \ref{L2}, \ref{L3} and the definition of $\rho$. (b) In the general case, the 
assertion follows from (a) and Part 1. \qed

\begin{cor}
$d(L(e_1), L(e_k))\ge 2\rho>0$. 
\end{cor}
\proof Note first that $x_1, x_{k+1}, e_1, e_k$ cannot be large. 
Since the (combinatorial) length of $\theta$ is at least $3$, 
by Lemma \ref{large exists} there exists a large cell $c$ in $\theta$ (every 
 edge $e_i$, $1<i<k$, will contain a large cell $c$). Then 
$$
e_1< c< e_k
$$   
and, by part 3 of the above lemma,
$$
L(e_1)\subset Bis(c)^-, L(e_k)\subset Bis(c)^+ 
$$
Moreover
$$
\min(d(L(e_1), Bis(c)), d(L(e_k), Bis(c)) )\ge \rho. 
$$
Corollary follows. \qed 

This concludes the proof of Theorem \ref{comb}. \qed

\begin{rem}
By using the arguments of the proof of Theorem \ref{comb}, one can also 
show that $\phi(G)$ is convex-cocompact.  
\end{rem}

\section{Arithmetic subgroups of $O(n,1)$}\label{sec:quat}

The goal of this section is to describe a ``quaternionic'' 
 construction of arithmetic subgroups in $O(n,1)$. For $n\ne 3, 7$, this construction covers all arithmetic subgroups. Our discussion follows \cite{Vinberg-Shvartsman}, we refer the reader to \cite{Li-Millson} 
 for the detailed proofs. 

We begin by reviewing quaterion algebras over number fields and ``hermitian vector spaces'' over such algebras.

Let $K$ be a field, $D$ be a central quaternion algebra over $K$. In other words, the algebra $D=D(a,b)$ has the basis $\{1, i, j, k\}$, subject to the relations:
$$
i^2=a\in K,\quad j^2=b\in K,\quad ij=-ji=k
$$
and so that $1$ generates the center of $D$. For instance, for $a=b=-1$ and $K=\R$, we get the algebra of Hamilton's quaternions ${\mathbf H}$.  Similarly, if $a=b=1$ then $D$ is naturally isomorphic to the algebra of $2\times 2$ matrices $End(\R^2)$.  One uses the notation
$$
\la= x1+ yi + zj +wk= x+ yi + zj +wk
$$
for the elements of $D$, with $x, y, z, w\in K$. An element of $D$ is {\em imaginary} if $x=0$. 
We will identify $K$ with the center of $D$: 
$$
K=K\cdot 1\subset D.
$$ 

One defines the {\em conjugation} on $D$ by 
$$
\la= x+ yi + zj +wk\mapsto \bar\la= x- yi - zj -wk. 
$$

Then $Tr(\la)= \la+ \bar\la$, $N(\la)=\la\bar\la$ are the {\em trace} and the {\em norm} on $D$. 
Clearly, both trace and the norm are elements of $K$. 
In the case when $D\cong End(\R^2)$, the trace is twice the matrix trace and the norm is the matrix determinant. Suppose that $K$ is a subfield of $\R$. 
We say that $\la\in D$ is {\em positive} (resp. {\em negative}) if it has positive (resp. negative) norm. 

\medskip
In what follows, we will assume that $K$ is a totally real number field and $D$ is a division algebra.

We will consider finite-dimensional ``vector spaces'' $V$ over $D$, i.e., finite-dimen\-sional right $D$-modules where we use the notation 
$$v\la=v\cdot \la \in V$$
for $v\in V, \la\in D$. Such a module is isomorphic to $D^n$ for some $n<\infty$, where $n$ is the dimension of $V$ as a $D$-module. Given $v \in  V$ define $<v>$ as the submodule in $V$ 
generated by $v$.

A {\em skew-hermitian form} 
on $V$ is a function $F(u, v)=\<u, v\>\in D$, $u, v\in V$, so that:
$$
F(u_1+u_2, v)= F(u_1,v)+ F(u_2,v),
$$
$$
F(u \la, v\mu)= \bar\la F(u,v) \mu, \quad F(u,v)= -\ol{F(v,u)}. $$
Similarly, $F$ is hermitian if 
$$
F(u \la, v\mu)= \bar\la F(u,v) \mu, \quad F(u,v)= \ol{F(v,u)}. $$
The form $F$ is {\em nondegenerate} if 
$$
 F(v, u)=0 \quad \forall u\in V \Rightarrow v=0.
$$
In coordinates:
$$
F(x,y)= \sum_{l,m} \bar{x}_l a_{lm} y_m, \quad a_{lm}= - \ol{a_{ml}}.
$$
In particular, the diagonal entries of the Gramm matrix of $F$ are imaginary. A {\em null-vector} is a vector with $F(v,v)= 0$, equivalently, $F(v,v)$ is not imaginary. We say that a vector $v$ is {\em regular} if it is not null. From now on, we fix $F$. 

Define $U(V, F)$, the group of {\em unitary} automorphisms of $V$, i.e., invertible endomorphisms which preserve $F$.

For a regular vector $v$, define the submodule $v^\perp$ in $V$:
$$
v^\perp=\{ u\in V: \<u, v\>=0\}.
$$
We will see below that
$$
V= <v> \oplus v^\perp
$$
and the restriction of the form $F$ to $v^\perp$ is again nondegenerate.

\medskip 
{\bf Orthogonal projection.} Suppose $\<v, v\>=a\ne 0$. Define
$$
Proj_v: V\to <v>, Proj_v(u)= v \cdot a^{-1} \<v, u\>.
$$
Then $Proj_v\in End(V)$, $Proj_v|_{<v>}=Id$, $Ker(Proj_v)=v^\perp$. In particular, for a vector $u\in V$,
$$
u'=u- Proj_v(u)\in v^\perp.
$$
Hence, $u=u'+u''$, with $u''=Proj_v(u)$.  It is now immediate that
$$
V= <v> \oplus v^\perp.
$$
Since $F$ is nondegenerate, the restriction $F|_{v^\perp}$ is also nondegenerate.

\medskip 
The existence of projections allows us to define Gramm-Schmidt  orthogonalization in $V$. In particular, $V$ has an orthogonal basis in which $F$ is diagonal.

\medskip 
{\bf Reflections.} Given a regular vector $v\in V$,  define the {\em reflection} $\si_v \in End(V)$ by
$$
\si_v(u):= u- 2 Proj_v(u).
$$
It is immediate that $\si_v|<v>=-Id$ and $\si_v| v^\perp = Id$. In particular, $\si_v$ is an involution  
in $U(V, F)$.

Observe that reflections $\si_{u}, \si_v$ commute iff either $<u>=<v>$ or $\<u, v\>=0$.

Proof of noncoherence of arithmetic lattices will use the following technical result:

\begin{lem}\label{L}
Let $V$ be 3-dimensional, $p_1, p_2$ be regular vectors which span a 2-dimensional submodule $P$ in $V$ so that the restriction of $F$ to $P$ is nondegenerate. Then there exist $u_1, u_2\in V$ so that:

1. $\<p_m, u_m\>=0, m=1, 2$.

2. $\<u_1, u_2\>=0$.
\end{lem}
\proof Since $F|_P$ is nondegenerate, it follows from the Gramm-Schmidt  orthogonalization that there exists
$v\in V$ so that $P= v^\perp$ and $V= P \oplus <v>$. In particular, $v$ is  a regular vector.

Orthogonalization implies that there exist vectors $u_1', u_2'\in P$ orthogonal to $p_1, p_2$ respectively. We will
find  vectors $u_1, u_2$ in the form
$$
u_m= u_m'+ v \cdot \la_m, \quad \la_m\in D, m=1, 2,
$$

It is immediate that $\< p_m, u_m\>=0$, $m=1, 2$. We have
$$
\< u_1, u_2\>= \< u_1', u_2'\> + \bar{\la}_1 \< v, u_2'\> + \< u_1', v\> \la_2 + \bar{\la}_1 \<v,v\> \la_2.
$$
Set  $\al:=\<v,v\>$, $\nu_1:=\< u_1', v\>, \nu_2:=  \< v, u_2'\>, \mu:= \< u_1', u_2'\> $. By scaling $u_1'$ if necessary, we get
$$
\nu_1+ \al\ne 0.
$$
We now set $\la_1=1, \la_2=\la$ (the unknown). Then the equation $\< u_1, u_2\>=0$ has the solution
$$
\la= - (\mu+ \nu_2)(\nu_1 +\al)^{-1}\in D. \qed
$$

\medskip 
Next, we now relate $(V, F)$ to hyperbolic geometry. Regarding $K$ as a subfield of $\R$, we define the completions $D_\R$ of $D$ and $V_\R$ of $V$. We require $D$ to be such that $D_\R \cong SL(2, \R)$, i.e., to have zero divisors. Hence,
at least one of the generators $i, j, k$ of $D$ is negative, i.e., has negative norm. We assume that this is $i$;
thus $i^2>0$. By interchanging $j$ and $k$ if necessary, we obtain that $j$ is also negative. By scaling $i, j$ by appropriate real numbers we get $i^2=j^2=1$. By abusing the terminology,
we retain the notation $i, j$ for these real multiples of the original generators of $D$.

Since $i^2=1$,
the right multiplication by $i$,
$$
I(v)= v\cdot i, v\in V,
$$
determines an involutive linear transformation of $V_\R$ (regarded as the real vector space). We obtain the eigenspace decomposition
$$
V_\R= V_+ \oplus V_-,$$
where $I|V_\pm = \pm Id$. (Note that $V_\pm \cap V=0$.) Let $J: V_\R\to V_\R$ be given by the right multiplication by $j$: This is again a linear automorphism. For $v\in V_+$ we have
$$
v j i= - vi j = -vj.
$$
Therefore, $J$ determines an isomorphism $V_+\to V_-$.

Our next goal is to analyze the subspace $V_+$. For $u, v\in V_+$ we have
$$
c=\< u, v\>= -i \<u, v\>= \<u i, v\>= \< u, v i\>= \<u, v\> i.
$$
Hence,
$$
-ic= c= ci.
$$
Such $c$ necessarily has the form $c= t(k-j), t\in \R$. Set $\al:= k-j$. Then $F|_{V_+}$ 
takes values in $\R\al$. In particular, $F|_{V_+}= \al \varphi$, where $\varphi$ is a real bilinear form.

Define $\be:=1+i$. Then
$$
\al i = \al, i\be=\be i= \be, \al^2=0, \be \bar \be=0, \al\be= \bar\be \al= 2\al,
$$
$$
k\be = \al, k \al= -\be.
$$

We now consider the case when $V$ is 1-dimensional, i.e., $V=D$.
Then the form $F$ is given by
$$
F(x,y)= \bar{x} a y, a\in D, a=- \bar{a}.
$$
Therefore, 
$$
D_+= \{ x \be + y \al: x, y\in \R\}.
$$
We now compute the form $\varphi$ so that $F|D_+= \al \varphi$. Note that the group
$$
SL_1(D)=\{g: N(g)=g \bar g=1\}$$
 acts on the space of traceless matrices
 $$
 D_0=\{ \la\in D_\R : Tr(\la)=\la+\bar\la=0\}$$
  by
 $$
 Ad_g(\la)= \bar{g} \la g= g^{-1} \la g.
 $$
 This action preserves the nondegenerate indefinite quadratic form $\la \bar \la$. Hence, this is a the orthogonal  action on $\R^{2,1}$ which has three nonzero orbit types. The relevant ones are positive ($N(\la)>0$) and negative ($N(\la)<0$) vectors in $D_0$. They are represented by $\la=k$ (in which case $N(k)=1$) and  $\la= i$ (in which case $N(i)=-1$).

 By changing the generator in $D_\R $, we replace $a$ (in the definition of $F$) with $\bar{g} a g$. Hence, our analysis of the form $\varphi$ reduces to two cases: $a=i, a=k$.

{\bf  Case 1:} $a=i$, $N(a)<0$. Then for $v=a \be + y \al \in D_+$, we have
 $$
F(v,v)= \< x \be + y \al, x \be + y \al\>= (x \bar\be + y \bar\al)i(x \be + y \al)=$$
$$
(x \bar\be - y\al)(x \be - y \al)= -2xy \al.
 $$
 Hence, $\varphi(v,v)= -2xy$, an indefinite form.

{\bf Case 2:} $a=k, N(a)>0$. Then for $v=a \be + y \al \in D_+$, we have
 $$
F(v,v)= \< x \be + y \al, x \be + y \al\>= (x \bar\be + y \bar\al)k(x \be + y \al)=$$
$$
(x \bar\be + y\bar\al)(x \al - y \be)= 2(x^2+y^2) \al.
 $$

 Therefore, in this case the form $\varphi$ is positive-definite.

 To summarize, if $N(a)<0$ then $\varphi$ is indefinite, while if $N(a)>0$ then the form $\varphi$ is positive-definite; in both cases, $\varphi$ is nondegenerate.

\medskip 
 We now consider the general case. Without loss of generality, we may assume that $F$ is diagonal,
 $\{v_1,...,v_n\}$ is an orthogonal basis in $V$. Clearly,
 $$
 V_\pm= \oplus_{l=1}^n <v_l>_\pm,
 $$
 where $<v_l>_\pm = v_l \cdot D_\pm$. We say that a vector $v\in V$ is {\em positive} (resp. {\em negative}) if $N(\<v,v\>)>0$ (resp. $N(\<v,v\>)<0$).

 Thus we obtain

 \begin{prop}
 The form $\varphi$ is always nondegenerate. It has signature $(p,q)$ iff $F$ has ``signature'' $(p,q)$, i.e.,
 some (equivalently, every) orthogonal basis $\{v_1,...,v_n\}$ of $V$  contains exactly $q$ negative vectors.

 In particular, $\varphi$ has signature $(p,1)$ iff one of the vectors of one (every)
 orthogonal basis of $V$ is negative and the rest of the vectors are positive.
 \end{prop}

 We assume from now on that $F$ has the signature $(p,1)$. Then $(V,F)$ determines the hyperbolic 
 $2p+1$-space $H(V)$ associated with the Lorentzian space $(V_+, \varphi)$. A vector $v\in V$ determines a geodesic $H(<v>)\subset H(V)$ iff $v$ is a negative vector. Similarly, $F|_{v^\perp}$ has hyperbolic signature iff $v$ is a positive vector; thus $H(v^\perp)\subset H(V)$ is a hyperbolic hyperplane.

For a subspace $W_+\subset V_+$ we define $W^-_+$ as the closed light-cone
$$
\{w\in W_+: \varphi(w)\le 0\}. 
$$

We obtain an embedding $U(V, F)\embed Isom(H(V))$ given by the action of $U(V, F)$ on the 
Lorenzian space $(V_+, \varphi)\cong \R^{p,1}$.

 \begin{lem}
 Suppose that $u, v$ are positive vectors in $V$. Then the hyperbolic hyperplanes $H(u^\perp), H(v^\perp)$ are orthogonal iff $u$ is orthogonal to $v$.
 \end{lem}
\proof Lemma immediately follows from the fact that the following are equivalent:

1. $H(u^\perp), H(v^\perp)$ are orthogonal or equal.

2. $\si_u$ commutes with $\si_v$.

3. $u$ and $v$ are orthogonal or generate the same submodule in $V$. \qed

\medskip
Given the hyperbolic space $H(V)$, we say that a hyperplane $L\subset H(V)$ is $K$-rational if $L=H(u^\perp)$ for some positive vector $u\in V$. A subspace in $H(V)$ is $K$-rational if it appears as intersection of $K$-rational hyperplanes. One observes (using orthogonalization) that a geodesic $\ga\subset H(V)$ is $K$-rational iff $\ga=H(<v>)$, where $v\in V$ is a negative vector.

We therefore obtain

\begin{cor}
Suppose that $\ga_1, \ga_2$ are distinct $K$-rational geodesics in the 
hyperbolic 5-space $\H^5= H(V)$, $dim(V)=3$. Then there exist $K$-rational hyperplanes $H_1, H_2$ containing $\ga_1, \ga_2$, so that $H_1$ is orthogonal to $H_2$.
\end{cor}
\proof Let $p_l\in V$ be such that $\ga_l=H(<p_l>), l=1, 2$. Then, according to Lemma \ref{L}, there exist $u_1, u_2\in V$ so that
$$
<p_l>\subset u_l^\perp, l=1, 2, \quad \<u_1, u_2\>=0.
$$
Since $F|_{<p_l>}$ is indefinite, it follows that $F|u_l^\perp$ is indefinite as well ($l=1,2$). Hence, $u_l^\perp$ determines a $K$-rational hyperbolic hyperplane $H_l=H(u_l^\perp)\subset \H^3$. Clearly,
$\ga_l\subset H_l$, $l=1, 2$. In view of the previous lemma, $H_1$ and $H_2$ are orthogonal provided that they actually intersect in $\H^5$.

Pick a generator $w$ of $W= u_1^\perp \cap u_2^\perp$. Since $u_1, u_2$ are positive and 
$u_1, u_2, w$ form an orthonormal basis in $V$, it follows that $w$ is negative. Therefore, $\ga=H(W)=H(<w>)$ is a geodesic in $\H^5$. Hence, $H_1, H_2$ intersect along the geodesic $\ga$  
at the right angle. \qed

\medskip
Let $\si_l:=\si_{u_l}, l=1, 2$ denote the reflections in the subspaces $U_l=u_l^\perp$. 
Since $\si_1, \si_2$ commute, their product is the involution $\tau$ whose fixed-point set is the 1-dimensional subspace $W=U_1\cap U_2$. We now assume that the geodesics $\ga_1, \ga_2$ are within positive distance from each other. Set $\ga:= H(W)$ (a $K$-rational geodesic in $\H^5$) and  set 
$L:= H(P)$ (a 3-dimensional $K$-rational hyperbolic subspace in $\H^5$).

\begin{lem}
$H$ and $\ga$, $L$ and $\tau(L)$ are within positive distance from each other.  
\end{lem}
\proof Observe that $U_l\ne P, l=1,2$ (since $u_l$ is not a multiple of $v$). In particular, 
$U_l\cap P=<p_l>$. Hence,
$$
U_{l+}^-\cap P_+^-=<p_l>_+^-. 
$$
If $w\in W_+^-\cap P_+^-$ then $w\in U_{l+}^-\cap P_+^-=<p_l>_+^-$, $l=1, 2$. However, 
$<p_1>_+^-\cap <p_2>_+^-=0$ since we assumed that $\ga_1, \ga_2$ are within positive distance 
from each other. Thus $w=0$. In particular, $L$ and $\ga$ are within positive distance from each 
other. 

Let $\rho\subset \H^5$ denote the geodesic segment with the end-points in $L, \ga$, which is orthogonal to both. Then $\rho\cup \tau(\rho)$ is a geodesic segment connecting $L$ and $\tau(L)$ and orthogonal to bot subspaces. Hence,  $L$ and $\tau(L)$ are within positive distance from each other. \qed

\medskip 
Suppose that $F$ is a skew-hermitian form on $V$. Given an embedding $\tau: K\to \R$ we define 
 the signature $sig_\tau(F)$ with respect to the subfield $\tau(F)\subset \R$: 
 Note that the notion of positivity and negativity in $V$ (and, hence, the signature) depends on the 
embedding  $\tau$. 

At last, we are ready to define the class of ``quaternionic'' arithmetic lattices  $G\subset O(n,1)$. Let $K$ be totally real, $D$ be a quaternion algebra over $K$ and $V$ be an $n+1$-dimensional module over $D$. We assume that $F$ is a nondegenerate hermitian form on $V$ satisfying the following:

1. $F$ has the signature $(n,1)$. 

2. For every nontrivial embedding $\tau: K\to \R$, the signature $sig_\tau(F)$ is $(n+1,0)$ 
(or $(0, n+1)$).  

Next, we need a notion of an ``integer'' automorphism of $V$. Let $O\subset D$ be an {\em order}, i.e., 
a lattice in $D$ regarded as a vector space over $K$.  An example of such order is given by 
$A^4\subset D$, where $A$ is the ring of integers of $K$. 

The order $O$ also determines the lattice $O^{n+1}\subset V=D^{n+1}$. We let $GL(V, O)$ denote the group of automorphisms of the $D$-module $V$ which preserve the lattice $O^{n+1}$. If we regard automorphisms of $V$ as ``matrices'' with coefficients in $D$, then the elements of $GL(V, O)$ are ``matrices'' with coefficients in $O$ which admit inverses with the same property. 

We now fix an order $O\subset D$. Then every subgroup $\Ga$ of $U(V, F)$ commensurable to the intersection $U(V, F)\cap GL(V, O)$ is {\em an arithmetic group of quaternionic type}. 
The embedding $U(V, F)\embed O(n,1)$ 
(induced by the identity embedding $K\embed \R$) realizes $\Ga$ as a lattice in $O(n,1)$. 

\begin{thm}
(See \cite{Li-Millson}.) 
Except for $n=3, 7$, every arithmetic subgroup of $O(n,1)$ appears as one of the groups $\Ga$ 
as above. In the case $n=7$, there is an extra class of arithmetic groups associated with octaves rather than quaternions. 
For $n=3$, there is yet another construction, also of quaternionic origin, which covers all arithemetic groups in this dimension, 
see \cite{Maclachlan-Reid} for the detailed description. 
\end{thm}

\section{Proof of noncoherence of arithmetic groups of quaternionic origin}\label{final}

Let $G_0\subset O(3,1)$ be an arithmetic lattice. 
According to our assumptions, there exists a finite index (torsion-free) 
subgroup $G_0'\subset G_0$ so that the manifold $M^3=\H^3/G_0'$ fibers over the circle. 
Let $F_0\triangleleft G_0'$ denote the normal surface subgroup corresponding to the fundamental  
group of the surface fiber in this fibration. Pick two nonconjugate maximal cyclic subgroups 
$G_{e_1}', G_{e_2}'$ in $G_0'$ so that $G_{e_i}'\cap F_0=\{1\}$. 

\begin{lem}
We can choose $G_{e_1}', G_{e_2}'$ so that the pair 
$\{G_{e_1}', G_{e_2}'\}$ is relatively separable in $G_0$. 
\end{lem}
\proof According to discussion in section \ref{RFRS}, 
there are at least two distinct virtual fibrations of the manifold 
$M^3$. Let $J_0\triangleleft G_0'$ be a normal surface subgroup corresponding to the second fibration. Then $J_0\cap F_0$ has infinite index in both $F_0, J_0$.  Therefore, take two nontrivial 
elements $t_1, t_2\in J_0$ so that:
$$
\<t_1\>\ne \<t_2\>, \quad \<t_1\>\cap F_0= \<t_2\>\cap F_0=\{1\}.  
$$
After passing to a finite-index subgroup in $G_0'$ if necessary, we can assume that 
the subgroups $G_{e_i}':=\<t_i\>, i=1,2$ are maximal cyclic subgroups of $G_0'$ and that they are not conjugate in $G_0'$. The set $\{\<t_1, t_2\>\}$ is relatively separable in $J_0$ by Proposition \ref{prop:sep} since, being a surface group, $J_0$ is LERF  \cite{Scott(1978)}. The subgroup $J_0$ is normal in $G_0'$  with $G_0'/J_0\cong \Z$. Therefore, relative separability of 
$\{\<t_1, t_2\>\}$ in $G_0'$  follows from Lemma \ref{easy}. \qed

\medskip 
We let $\ga_i, i=1,2$ denote the invariant geodesic 
of $G_{e_i}'$. Let $R_0:=d(\ga_1, \ga_2)$. We will assume, as in the beginning of Section 
\ref{comb} that $R_0$ is the distance between the projections of $\ga_1, \ga_2$ to the manifold $M^3$.  

Since $G_{e_1}', G_{e_2}'$ are separable in $G_0'$, 
given a number $R_1$, there exists a finite-index subgroup $G_0''\subset G_0'$ so that:
 
For each $\ga=\ga_i, i=1, 2$ and geodesics 
$\be_1, \be_2\in G_0''\cdot (\ga_1 \cup \ga_2)$ 
so that $d(\be_j, \ga)= R_0$ ($j=1, 2$), it follows that 
$$
d(proj_{\ga}(\be_1), proj_{\ga}(\be_2))\ge R_1.
$$

Set $G_{e_i}'':=G_{e_1}'\cap G_0'', i=1, 2$. Without loss of generality, we may assume 
that $G_{e_1}'', G_{e_2}''$ generate a free subgroup $H$ of $F_0'':=F_0\cap G_0''$. The subgroup 
$H$ is separable in $F_0''$; since $G_0''/F_0''\cong \Z$, the subgroup $H$ is also separable in $G_0''$. 
Therefore, in view of Lemma \ref{easy}, the pair $\{G_{e_1}'', G_{e_2}''\}$ is weakly separable in $G_0''$. 
Therefore, given a number $R$, one can find a finite-index subgroup $G_0'''\subset G_0''$ so that:
 
For all distinct geodesics $\be, \ga$ in the $G_0'''$-orbit of $\ga_1, \ga_2$, the distance 
$d(\be, \ga)$ is at least $R$ unless there exists $g\in G_0'''$ 
which carries $\be\cup \ga$ to $\ga_1\cup \ga_2$. 

Therefore, the Assumption \ref{as1} (from section \ref{comb}) is satisfied by the group $G_1:=G_0'''$ and 
its subgroups $G_{\eps_i}:=G_{e_i}\cap G_1$, $i=1,2$, with respect to the numbers $R_1$ and $R$. 

Our next goal is to define $R_1$ and $R_4$. We let $S_i, i=1,2$, be orthogonal 3-dimensional 
$K$-rational subspaces in $\H^5$ intersecting $L_1$ along the geodesics $\ga_1, \ga_2$; 
let $L_2:=\si_s(L_1), L_4:=\si_1(L_1)$. Since $S_1, S_2$ are $K$-rational, the 
involutions $\si_1, \si_2$ belong to the commensurator of the lattice $\Ga$. Since the group generated by $\si_1, \si_2$ is finite, without loss of generality we may assume that these involutions normalize $\Ga$ 
(otherwise, we first pass to a finite-index subgroup in $\Ga$). 

Then $r>0$ is the distance $d(L_2, L_4)$. Let $\al_1, \al_2$ denote the angles $\angle(L_1, L_4)$, $\angle(L_1, L_2)$ and $\al:=\min(\al_1, \al_2)$. 

In view of Lemma \ref{L*}, we will use $R\ge R_*$ so that for every $g\in G_1$,
$$
d(g(L_i), L_j)\ge r, i, j\in \{2,4\}, g(L_i)\ne L_j. 
$$

Lastly, we set $R_1:=R_1(\al, R_0, r)$ and $R_4:=R_4(\al, R_0, r)$. We then will use $R:=R_4$.

As in section \ref{sec:comb}, we define a quadrilateral $Q$ of groups with vertex groups isomorphic to $G_1$, so that these isomorphisms send $G_{\eps_1}, G_{\eps_4}$  to edge groups. Let 
$G:=\pi_1(Q)$. Using the involutions $\si_1, \si_2$ as in section \ref{sec:comb}, we construct a 
discrete and faithful representation $\phi: G\to O(n,1)$. Since $G_1\subset \Ga$ and $\si_1, \si_2$ normalize $\Ga$, the image of this representation is contained in $\Ga$. (This provides yet another proof of discreteness of $\phi(G)$, however we still have to use Combination Theorem \ref{comb} in order to conclude that $\phi$ is faithful.  

In order to show incoherence of $\Ga$ it suffices to prove:

\begin{lem}
The group $G$ is noncoherent. 
\end{lem}
\proof Let $\eps_1,...,\eps_4$ denote the edges of $Q$ and $G_{\eps_i}, i=1,...,4$ denote the corresponding edge groups. Recall that $F_1\triangleleft G_1$ is a normal surface subgroup. 
Let $F_i$ denote the normal surface subgroups of $G_i, i=2,3,4$, which are the images of $F_1$  under the isomorphisms $G_1\to G_i$. Then 
$$
F_i\cap G_{\eps_{i-1}}=F_i\cap G_{\eps_i}=\{1\}, 
$$
here and in what follows $i$ is taken mod 4. Let $F$ denote the subgroup of $G$ generated by $F_1,...,F_4$. Clearly, this group is finitely-generated. We will show that $F$ is not finitely presented 
by proving that $F\cong F_+ *_N F_-$, where $F_\pm$ are finitely-generated (actually, finitely-presented)  and $N$ is a free group of infinite rank.  

We first describe $G$ as an amalgamated free product: We cut the quadrilateral $Q$ in half so that one half contains the vertices $1, 2$, while the other half contains the vertices $3,4$. 
Accordingly, set 
$$
G_-:=\<G_1, G_2\>\cong G_1*_{G_{\eps_2}} G_2,\quad G_+:=\<G_3, G_4\>\cong G_3*_{G_{\eps_4}} G_4, 
$$
$$
E:= \<G_{\eps_1}, G_{\eps_3}\>\cong G_{\eps_1}* G_{\eps_3}\cong \Z * \Z. 
$$
Then 
$$
G\cong G_- *_E G_+. 
$$
Similarly, we set
$$
F_+:=F\cap G_+\cong F_1* F_2; \quad F_-:=F\cap G_-\cong F_3* F_4. 
$$
Since $G_1$ is generated by $t_2$ and $F_1$, the group $G_2$ is generated by $t_2$ and $F_2$. 
Hence, $F_-$ is normal in $G_-$ and $G_-/F_-\cong \Z$. Moreover, $F_-$ has trivial intersection with 
$G_{\eps_1}$ (since $F_1$ does). It is immediate that $N:=F_-\cap E$ is an infinite index nontrivial subgroup of $E$. Since $E$ is free of rank 2, it follows that $N$ is a free group of infinite rank. 
Clearly, $N=F_+\cap E$ and we obtain
$$
N\cong F_+ *_N F_-. 
$$

Therefore, $F$ is finitely generated 
and infinitely presented since (by considering the Meyer-Vietoris sequence associated with the amalgam 
$N\cong F_+ *_N F_-$) $H_2(F, \Z)$ has infinite rank. Thus, $G$ is noncoherent. \qed 

 \section{Complex-hyperbolic and quaternionic lattices}\label{sec:ch}

It is an important open problem in theory of lattices in rank 1 Lie groups $O(n,1)$ and $SU(n,1)$ if 
a lattice has positive virtual first Betti number, i.e., contains a finite-index subgroup with infinite abelianization. In this section we relate this problem to noncoherence in the case of $SU(n,1)$. It was proven by D.~Kazhdan \cite{Kazhdan(1977)} 
(see also \cite{Wallach(1984)}) that arithmetic lattices of the {\em simplest type} (or, first type) 
in $SU(n,1)$ admit finite index (congruence) subgroups with infinite abelianization. 
Certain classes of non-arithmetic lattices in $SU(2,1)$ (the ones violating integrality condition for arithmetic groups) are proven to have positive virtual first 
Betti number by S.-K. Yeung \cite{Yeung(2004)}. 

On the other hand J.~Rogawski \cite{Rogawski} proved that for arithmetic lattices $\Ga$ in $SU(2,1)$ 
of {\em second type} (associated with division algebras), 
every congruence-subgroup $\Ga'\subset \Ga$ has finite abelianization. 
It is unknown if non-congruence subgroups in such lattices (if they exist at all!) can have infinite abelianization. 

Below is the description of arithmetic lattices of the simplest type in $SU(n,1)$ following 
\cite{McReynolds} and  \cite{Stover}.  
Let $K$ be a totally real number field;  %of degree $k$ over $\Q$ 
take a totally imaginary quadratic extension $L/K$ and let $\O_L$ be the ring of integers 
of $L$. Let $\si_1, \si_2,..., \si_k: L\to \C$ be the embeddings.  
%so that 
%$$\si_j|K= \hat\si_j|K, j=1,...,k. $$ 
Next, take a hermitian quadratic form in $n+1$ variables 
$$
\varphi(z, \bar{z})=\sum_{p,q=1}^{n+1} a_{pq} z_p \bar{z}_q
$$
with coefficients in $L$. We require $\varphi^{\si_1}, \varphi^{\bar\si_1}$ to have signature $(n,1)$  
and require the forms $\varphi^{\si_j}$ to have signature $(n+1,0)$ for the rest of the embeddings 
$\si_j$. Let $SU(\varphi)$ denote the group of special unitary automorphisms of the form $\varphi$ on 
$L^{n+1}$. The embedding $\si_1$ defines a homomorphism $SU(\varphi)\to SU(n,1)$ with relatively compact kernel. We will identify $L$ with $\si_1(L)$, so $\si_1=id$. 

\begin{defn}
A subgroup $\Ga$ of $SU(n,1)$ is said to be an arithmetic lattice of the {\em simplest type} if it is commensurable to $SU(\varphi, \O_L)= SU(\varphi)\cap SL(n+1, \O_L)$. 
\end{defn}

By diagonalizing the form $\varphi$, we see that $L^{n+1}$ contains a 3-dimensio\-nal 
subspace $L^3$ so that the restriction of $\varphi$ to $L^3$ is a form of the signature $(2,1)$. 
Therefore, an arithmetic lattice $\Ga\subset SU(n,1)$ of the simplest type intersects $SU(\varphi|L^3)$ 
along a lattice $\Ga'$ of the simplest type (regarded as a subgroup of $SU(2,1)$). If $\Ga'$ is non-coherent, so is $\Ga$. Therefore, we restrict our discussion to the case of isometries of the complex-hyperbolic plane $\C \H^2$.
 
Suppose that $\Ga\subset SU(2,1)$ is a torsion-free uniform lattice with infinite abelianization. 
Therefore, $b_1(M)>0$, where $M=\C \H^2/\Ga$, 
Since $M$ is K\"ahler, its Betti numbers are even; therefore, 
there exists an epimorphism $\psi: \Ga\to \Z^2$. There are two cases to consider:

\medskip 
Case 1. $Ker(\psi)$ is not finitely generated. Then, according to \cite{Delzant(2006)}, 
there exists a holomorphic fibration $M\to R$ with connected fibers, where $R$ is a hyperbolic Riemann surface-orbifold. It was proven in \cite{Kapovich(1998a)} that the kernel $K$ of the homomorphism 
$\Ga\to \pi_1(R)$ is finitely-generated but not finitely presented. Hence, $\Ga$ is non-coherent in this case. 

\begin{rem}
Jonathan Hillman \cite{Hillman(letter)} suggested an alternative proof that $K$ is not finitely presented. Namely, 
if $K$ is of type $FP_2$ (e.g., is finitely presented) then it is a $PD(2)$-group 
(see Theorem \ref{hillman}) and, hence,  a surface group. 
It was proven by Hillman in \cite{Hillman(2000)} that the  holomorphic fibration 
$M\to R$ has no singular fibers. Such fibrations cannot exist due to a result of K.~Liu \cite{Liu(1996)}. 
Thus, $K$ is not finitely presented. 
\end{rem}

\medskip 
Case 2. $F=Ker(\psi)$ is finitely generated. If $\Ga$ were coherent, $F$ would be also finitely-presented. It is proven by Jonathan Hillman (Theorem \ref{hillman}) 
that $F$ has to be a surface group. We obtain the associated homomorphism 
$\eta:\Z^2\to Out(F)$ (the mapping class group of a surface). Since $\Ga$ contains no rank 2 abelian subgroups, 
$\eta$ is injective. Rank 2 abelian subgroups of the mapping class group have to contain nontrivial 
reducible elements \cite{Birman-Lubotzky-McCarthy}. 
Let $\ga\in \Z^2\setminus \{1\}$ be such that $\eta(\ga)$ is 
a reducible  element of the mapping class group. Hence, $\eta(\ga)$ fixes a conjugacy class of 
some $\al\in F\setminus \{1\}$. It follows that $\Ga$ contains $\Z^2$ (generated by a lift of $\ga$ to 
$\Ga$ and by $\al$). Contradiction. 

We thus obtain:

\begin{thm}
Suppose that $\Ga\subset SU(2,1)$ is a cocompact arithmetic group with infinite abelianization.  
Then $\Ga$ is non-coherent. 
\end{thm}

\begin{cor}
Suppose that $\Ga\subset SU(n,1)$ is a cocompact arithmetic group of the simplest type, where 
$n\ge 2$. Then $\Ga$ is non-coherent. 
\end{cor}

We now consider quaternionic-hyperbolic lattices. Recall that all lattices in ${\mathbf H} \H^n$, $n\ge 2$, are arithmetic accoriding to \cite{Corlette, Gromov-Schoen}. On the other hand, 
$Isom({\mathbf H} \H^1)\cong Isom(\H^4)$ and, hence, this group 
contains nonarithmetic lattices as well. 

\begin{prop}
Every arithmetic lattice in ${\mathbf H} \H^n$ is non-coherent. 
\end{prop}
\proof According to \cite{Platonov-Rapinchuk}, 
all arithmetic lattices in  $Isom({\mathbf H} \H^n)\cong Sp(n,1)$ have the following form. 

Let $K\subset \R$ be a totally-real number field, $D$ be a central quaternion 
algebra over $K$, $V$ be an $n+1$-dimensional right $D$-module, $F$ be a hermitian bilinear form on $V$ (see section \ref{sec:quat}). Choose a basis where $F$ is diagonal:
$$
F(x,y)=\sum_{m=1}^{n+1} \bar{x}_m a_{m} y_m,  
$$
$a_{m}=\bar{a}_{m}$. Then the signature of $F$ is $(p,q)$ if (after permuting the coordinates) 
$a_m>0, m=1,...,p$ and $a_{m}<0, m=p+1,...,n+1=p+q$. Let 
$U(V,F)$ be the group of unitary transformations of $(V,F)$. 

Given an embedding $\si: K\to \R$, we define a new form $F^\si$ 
by applying $\si$ to the coefficients of $F$. We now require $F, D$ and $K$ to be such 
that:

1.  $F$ has signature $(n,1)$ and $F^\si$ is definite for all embeddings $\si$ different from the identity. 

2. The completions of $D$ with respect to all the embeddings $\si: K\to \R$ are isomorphic to Hamilton's quaternions ${\mathbf H}$ (i.e., are division algebras).  

In particular, the embedding 
$D\to {\mathbf H}$, induced by the identity embedding $K\embed \R$, 
 gives rise to a homomorphism $\eta: U(V,F)\to Sp(n,1)=Isom({\mathbf H}\H^n)$. 

Let $O$ be an order in $D$ and set $\Ga_{V,O}:= U(V, F)\cap SL(V,O)$. 
 Lastly,  a group commensurable to $\eta(\Ga_{V,O}) \subset Sp(n,1)$ 
is called an {\em arithmetic} lattice in $Sp(n,1)$. Note that the kernel of the homomorphism
$$
\eta: \Ga_{V,O}\to Sp(n,1)
$$
is finite. Hence, very arithmetic lattice in $Sp(n,1)$ is abstractly commensurable to $\Ga_{V,O}$ for some choice of $K, D$ and $O$.

By restricting the form $F$ to the 2-dimensional submodule $W$ in $V$ spanned by the first and last basis vectors, we obtain a hermitian form of signature $(1,1)$. Therefore, every arithmetic lattice $\Ga$ in 
$Sp(n,1)$ will contain a subgroup commensurable to $\eta(\Ga_{W,O})$. The latter is an arithmetic lattice in $\H^4$ and, hence, is noncoherent according to \cite{KPV} and \cite{Agol(2008)}. Thus, $\Ga$ is incoherent as well. \qed  

The same argument applies to lattices $\Ga$ in the isometry group of the hyperbolic  plane over Cayley octaves $Isom({\mathbf O} \H^2)$, as every such lattice is arithmetic and contains an arithmetic 
sublattice  $\Ga'\subset Isom({\mathbf O} \H^1)\cong Isom(\H^8)$. 
(I owe this remark to Andrei Rapinchuk.) Since $\Ga'$ is noncoherent, so is $\Ga$. 

\bibliographystyle{siam}
\bibliography{lit}

\end{document}